\documentclass[graybox]{svmult}
\usepackage{amsmath,amssymb,bm}
\usepackage{mathptmx}       
\usepackage{helvet}         
\usepackage{courier}        
\usepackage{type1cm}        
\usepackage{makeidx}         
\usepackage{graphicx}        
\usepackage{multicol}        
\usepackage[bottom]{footmisc}
\makeindex
\newcommand{\jump}[1]{\left[\hspace{-0.025in}\left[#1\right]\hspace{-0.025in}\right]}
\newcommand{\avg}[1]{\left\{\hspace{-0.045in}\left\{#1\right\}\hspace{-0.045in}\right\}}

\newcommand{\figref}[1]{Fig.~\ref{#1}}

\newcommand{\circover}[1]{\mathaccent'27{#1}}

\newcommand{\semin}[1]{|#1|_{H^2(\Omega, {\mathcal{T}_h})}}
\def\bbR{\mathbb{R}}
\def\SS{H_{0}^{2}(\Omega)}
\def\L2{L_{2}(\Omega)}
\def\vecN{\mathbf{n}}
\def\Eh{\mathcal{E}_h}
\def\EhG{\mathcal{E}_{h,\Gamma}}
\def\Sz1{H_{0}^{1}(\Omega)}
\def\Vh{V_h}

\def\TEle{T}
\def\D{\Omega}
\def\p{\partial}
\def\BdyD{\p \D}
\def\VD{V_{h,D}}
\def\VC{V_{h,C}}

\def\cT{\mathcal{T}}
\def\Th{\cT_h}
\def\PCS{B_{\hspace{-1.5pt}B\!D\!D\hspace{-.5pt}C}}
\def\vD{v_{\scriptscriptstyle D}}
\def\vC{v_{\scriptscriptstyle C}}
\def\vCz{{v}_{\scriptscriptstyle C,\D\setminus\Gamma}}
\def\vCsc{{v}_{\scriptscriptstyle C,\Gamma}}

\def\ExtH{\mathcal{H}}
\def\Hcc{\mathcal{H}_{\mathcal{C}}}
\def\Hc{\mathcal{H}_{0}}
\def\Hcz{\mathaccent'27{\mathcal{H}}}

\def\ahC{a_{h}^C}
\def\AvgP{P_{\Gamma}}

\def\bbE{\mathbb{E}}
\def\vz{\circover{v}}
\def\AhG0{A_{h,C,\D\backslash\Gamma}}
\def\IhG0{I_{h,C,\D\backslash\Gamma}}

\def\PCS{B_{\hspace{-1.5pt}B\!D\!D\hspace{-.5pt}C}}
\def\VcG{V_{h,C}(\Gamma)}

\def\d{\displaystyle}
\def\ds{\displaystyle}
\begin{document}
\title*{A Balancing Domain Decomposition by Constraints Preconditioner for
 a $C^0$ Interior Penalty Method}
\titlerunning{A BDDC preconditioner for $C^0$ Interior Penalty Methods}
\author{Susanne C. Brenner, Eun-Hee Park, Li-Yeng Sung, and Kening Wang}
\institute{Susanne C. Brenner \at Department of Mathematics and Center for Computation and Technology, Louisiana State University, Baton Rouge, LA 70803, USA, \email{brenner@math.lsu.edu}\\ 
Eun-Hee Park \at School of General Studies, Kangwon National University, Samcheok, Gangwon 25913, Republic of Korea, \email{ eh.park@kangwon.ac.kr}\\ 
Li-Yeng Sung \at Department of Mathematics and Center for Computation and Technology, Louisiana State University, Baton Rouge, LA 70803, USA, \email{sung@math.lsu.edu}\\
Kening Wang \at Department of Mathematics and Statistics, University of North Florida, Jacksonville, FL 32224, USA, \email{ kening.wang@unf.edu}}
\maketitle
\abstract*{We develop a nonoverlapping domain decomposition preconditioner for the $C^0$ interior penalty method, a discontinuous Galerkin method, for the biharmonic problem. The preconditioner is based on balancing domain decomposition by constraints (BDDC). We prove that the condition number of the preconditioned system is bounded by $C (1+\ln (H/h))^2$, where $h$ is the mesh size of the triangulation, $H$ is the typical diameter of subdomains, and the positive constant $C$ is independent of $h$ and $H$. Numerical experiments are also represented to corroborate the theoretical result.}
%
\section{Introduction}\label{Sec:intro}
\noindent Consider the following weak formulation of a fourth order problem on a bounded polygonal domain $\D$ in $\bbR^{2}$:
\par\noindent
 Find $u \in \SS$ such that
\begin{equation}\label{eq:model-weakform}
  \int_\D \nabla^2 u : \nabla^2 v\,dx = \int_\D fv\,dx \qquad \forall\, v \in \SS,
\end{equation}
 where $f \in \L2$, and $\nabla^2 v : \nabla^2 w = \sum_{i,j=1}^{2} (\p^2v / \p x_i \p x_j)(\p^2w / \p x_i \p x_j)$ is the inner product of the Hessian matrices of $v$ and $w$.
\par
For simplicity, let $\Th$ be a quasi-uniform triangulation of $\D$ consisting of rectangles and take $\Vh \subset \Sz1$ to be the $Q_2$ Lagrange finite element space associated with $\Th$. Then the model problem \eqref{eq:model-weakform} can be discretized by the following $C^0$ interior penalty Galerkin method
 \cite{EGHLMT:2002:C0, BS:2005:C0IP}:
\par\noindent
 Find $u_{h}\in \Vh$ such that
\begin{equation*}\label{eq:SIPG}
  a_{h}(u_{h},v) = \int_\D fv\,dx \qquad v \in \Vh,
\end{equation*}
 where
\begin{align*}\label{eq:bilinear_form}
   a_h(v,w) & = \sum_{D \in \Th} \int_{T} \nabla^2 v : \nabla^2 w \, dx +
  \sum_{e\in \Eh} \frac{\eta}{|e|} \int_{e} \jump{\frac{\p v}{\p \vecN}} \jump{\frac{\p w}{\p \vecN}} \, ds
  \\
  & \hspace{40pt} +  \sum_{e\in \Eh} \int_{e} \left( \avg{\frac{\p^2v}{\p \vecN^2}}\jump{\frac{\p w}{\p \vecN}} +  \avg{\frac{\p^2w}{\p \vecN^2}} \jump{\frac{\p v}{\p \vecN}} \right)\, ds. \nonumber
\end{align*}
 Here $\eta$ is a positive penalty parameter, $\Eh$ is the set of edges of
 $\Th$, and $|e|$ is the length of the edge $e$. The jump $\jump{\cdot}$ and the average $\avg{\cdot}$ are defined as follows.
\begin{figure}
  \centering
    \includegraphics[height=2.5cm]{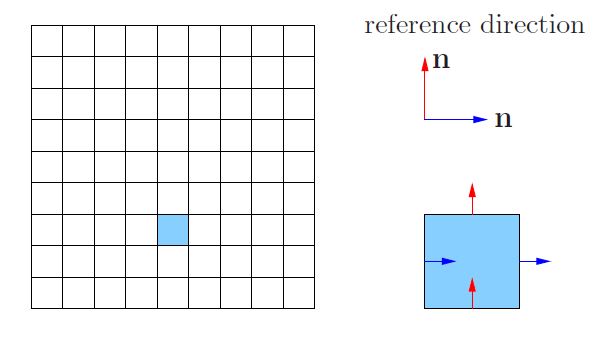}
    \caption{(a) A triangulation of $\D$. (b) A reference direction of normal vectors on the edges of $\TEle \in \Th$.}
    \label{fig:Th:ref_direction_n}
\end{figure}
\par
\noindent Let $\vecN_{e}$ be the unit normal chosen according to a reference direction shown in \figref{fig:Th:ref_direction_n}. If $e$ is an interior edge of $\Th$ shared by two elements $D_-$ and $D_+$, we define on $e$,
$$ \jump{\frac{\p v}{\p \vecN}} = \frac{\p v_{+}}{\p \vecN_{e}}-\frac{\p v_{-}}{\p \vecN_{e}} \quad \text{and} \quad
 \avg{\frac{\p^2 v}{\p \vecN^2}} = \frac{1}{2}\left(\frac{\p^2 v_{+}}{\p \vecN_{e}^2} + \frac{\p^2 v_{-}}{\p \vecN_{e}^2}\right),$$
\noindent where $v_{\pm} = v|_{D_{\pm}}$. On an edge of $\Th$ along $\BdyD$, we define
$$\jump{\frac{\p v}{\p \vecN}} = \pm \frac{\p v}{\p \vecN_e} \quad \text{and} \quad \avg{\frac{\p^2 v}{\p \vecN^2}} = \frac{\p^2 v}{\p \vecN_e^2},$$
in which the negative sign is chosen if $\vecN_e$ points towards the outside of $\D$, and the positive sign otherwise.
\par
It is noted that for $\eta>0$ sufficiently large (Lemma~6 in~\cite{BS:2005:C0IP}), there exist positive constants $C_1$ and $C_2$ independent of $h$ such that
\begin{equation*}\label{eq:seminorm_equiv}
   C_1 a_h(v,v) \leq \semin{v}^2 \leq C_2 a_h(v,v) \quad \forall v \in \Vh,
\end{equation*}
where 
\begin{equation*}
  \semin{v}^2 = \sum_{D \in \Th} | v |^{2}_{H^2(D)} + \sum_{e \in \Eh} \frac{1}{|e|}\left\|\jump{\frac{\p v}{\p \vecN}}\right\|^2_{L_2(e)}.
\end{equation*}
\par
Compared with classical finite element methods for fourth order problems, $C^0$ interior penalty methods have many advantages \cite{BS:2005:C0IP, BW:2012:BPSforC0, EGHLMT:2002:C0}. 
However, due to the nature of fourth order problems, the condition number of the discrete problem resulting from $C^0$ interior penalty methods grows at the rate of $h^{-4}$ \cite{LW:2007:DDM}. Thus a good preconditioner is essential for solving the discrete problem efficiently and accurately. In this paper, we develop a nonoverlapping domain decomposition preconditioner for $C^0$ interior penalty methods that is based on the balancing domain decomposition by constraints (BDDC) approach \cite{Dohrmann:2003:BDDC, BS:2007:BDDCFETI, BPS:2017:SIPGBDDC}.
\par
The rest of the paper is organized as follows. In
 Section~\ref{Sec:Precond_B1} we introduce the subspace
 decomposition.  We then design a BDDC preconditioner
 for the reduced problem in Section~\ref{Sec:Precond_B2}, followed by  condition number estimates in Section~\ref{Sec:CondNoEst}. Finally, we report numerical results in Section~\ref{Sec:NumericalResults} that illustrate the performance of the proposed preconditioner and
 corroborate the theoretical estimates.
\section{A Subspace Decomposition}\label{Sec:Precond_B1}
\noindent
We begin with a nonoverlapping domain decomposition of $\D$ consisting of rectangular (open) subdomains $\D_1, \D_2, \cdots, \D_J$ aligned with $\Th$ such that $\p \D_j \bigcap \p \D_\ell = \emptyset,$ a vertex, or an edge, if $j \neq \ell$.
\par
 We assume the subdomains are shape regular and denote the typical diameter of the subdomains by $H$. Let $\Gamma = \left( \bigcup_{j=1}^J \p\D_j \right) \backslash \p\D$ be the interface of the subdomains, and $\EhG$ be the subset of $\Eh$ containing the edges on $\Gamma$. 
%
\par
\begin{figure}
  \centering
    \includegraphics[height=2.5cm]{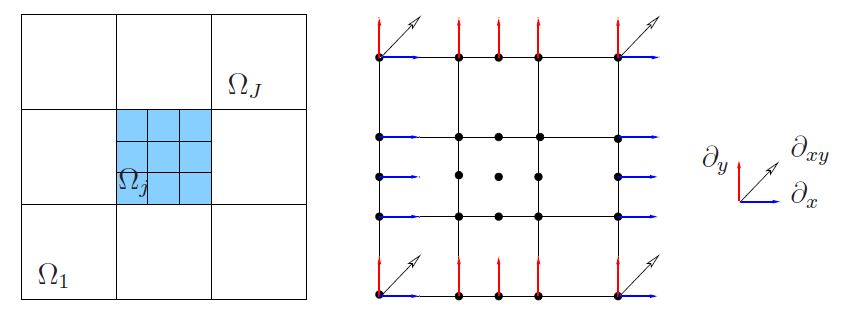}
\caption{(a) A nonoverlapping decomposition of $\D$ into $\D_1, \cdots, \D_J$ and a triangulation of the subdomain $\D_j$. (b) Dofs of $\Vh|_{\D_j}$. (c) Reference directions for the first order and mixed derivatives. }\label{fig:DD:dofsOfDj}
 \end{figure}

Since the condition that the normal derivative of $v$ vanishes on $\Gamma$ is implicit in terms of the standard degrees of freedom (dofs) of the $Q_2$ finite element, it is more convenient to use the modified $Q_2$ finite element space (\figref{fig:DD:dofsOfDj}) as $\Vh$. Details of the modified $Q_2$ finite element space can be found in \cite{BW:2012:BPSforC0}.
\par
First of all, we decompose $V_h$ into two subspaces
\begin{equation*}\label{eq:XhDecomposition}
  V_h = \VC \oplus \VD,
\end{equation*}
where
\begin{equation*}
   \VC = \bigg\{ v \in V_h : \jump{\frac{\p v}{\p \vecN}}=0 \text{ on the edges in $\Eh$ that are subsets of } \bigcup_{j=1}^J \p\D_j \bigg \} 
\end{equation*}  
and
\begin{align*} 
   \VD &= \bigg\{ v \in V_h : \avg{\frac{\p v}{\p \vecN}} = 0 \text{ on edges in } \EhG, \text{ and } \nonumber\\
  &\hspace{3.5cm} v \text{ vanishes at all interior nodes of each subdomain} \bigg \}.
\end{align*}
\par
Let $A_h : \Vh \rightarrow \Vh'$ be the symmetric positive definite (SPD) operator defined by
\begin{equation*}
 \langle A_h v, w \rangle = a_h (v, w) \qquad \forall\, v, w \in \Vh,
\end{equation*}
where $\langle \cdot, \cdot \rangle$ is the canonical bilinear form between a vector space and its dual. Similarly, we define $A_{h,C} : \VC \rightarrow \VC'$ and $A_{h,D} : \VD \rightarrow \VD'$ by
\begin{equation*}
 \langle A_{h,C} v, w \rangle = a_h (v, w) \quad \forall\, v, w \in \VC  \quad \text{and} \quad
 \langle A_{h,D} v, w \rangle = a_h (v, w) \quad \forall\, v, w \in \VD. 
\end{equation*}
\par
Then we have the following lemma.
\begin{lemma}\label{lm:operator_decomp}
 For any $v\in \Vh$, there is a unique decomposition $v = \vC + \vD$, where $\vC \in \VC$ and $\vD \in \VD$. In addition, it holds that
 \begin{equation*}
  \langle A_h v, v \rangle \approx \langle A_{h,C} v_C, v_C \rangle + \langle A_{h,D} v_D, v_D \rangle \quad \forall\, v\in \Vh.
 \end{equation*}
\end{lemma}
\par
\begin{remark}\label{rem:implementation_AhCAhD}
 Since the subspace $\VD$ only contains dofs on the boundary of subdomains, the size of the matrix $A_{h,D}$ is of order $J/h$. We can implement the solve $A_{h,D}^{-1}$ directly. Therefore, it is crucial to have an efficient preconditioner for $A_{h,C}$.
\end{remark}
\par
 Because functions in $\VC$ have continuous normal derivatives on the edges in $\EhG$ and vanishing normal derivatives on $\p\D$, it is easy to observe that
 \begin{equation*}\label{eq:localAjsum}
  a_h(v, w) = \sum_{j=1}^J a_{h,j} (v_j, w_j) \qquad \forall \, v, w \in \VC,
 \end{equation*}
 where $v_j = v\big|_{\D_j}, w_j = w\big|_{\D_j}$, and $a_{h,j}(\cdot, \cdot)$ is the analog of $a_h(\cdot, \cdot)$ defined on elements and interior edges of $\D_j$. Note that $a_{h,j}(\cdot, \cdot)$ is a localized bilinear form.
\par
Next we define
 \begin{align*}
  \VC (\D\backslash\Gamma) &= \left\{ v\in\VC : v \,\,\,\mbox{has its vanishing derivatives up to order 1 on} \,\,\,\Gamma \right\} \\
  \VC (\Gamma) &= \left\{ v\in\VC : a_h (v, w)=0 , \,\,\,\forall\, w\in\VC (\D\backslash\Gamma) \right\}. 
 \end{align*}
\par
Functions in $\VC (\Gamma)$ are referred to as discrete biharmonic functions. They are uniquely determined by the dofs associated with $\Gamma$.
\par
 For any $v_C \in \VC$, there is a unique decomposition $v_C = v_{C,\D\backslash\Gamma} + v_{C,\Gamma}$, where $v_{C,\D\backslash\Gamma} \in \VC (\D\backslash\Gamma)$ and $v_{C,\Gamma} \in \VC (\Gamma)$. Furthermore, let $A_{h,C,\D\backslash\Gamma} : \VC (\D\backslash\Gamma) \rightarrow \VC (\D\backslash\Gamma)'$ and $S_h : \VC (\Gamma) \rightarrow \VC (\Gamma)'$ be SPD operators defined by
\begin{align*}
 \langle A_{h,C,\D\backslash\Gamma} v, w \rangle &= a_h (v, w) \quad \forall \, v, w \in \VC (\D\backslash\Gamma), \\
 \langle S_h v, w \rangle &= a_h (v, w) \quad \forall \, v, w \in \VC (\Gamma), 
\end{align*}
then it holds that for all $\vC \in \VC$ with $\vC = \vCz + \vCsc$,
$$ \langle A_{h,C} v_C, v_C \rangle = \langle A_{h,C,\D\backslash\Gamma} v_{C, \D\backslash\Gamma}, v_{C, \D\backslash\Gamma} \rangle + \langle S_h v_{C, \Gamma}, v_{C, \Gamma} \rangle.$$
\par
\begin{remark}
{\rm  It is noted that $A_{h,C,\D\backslash\Gamma}^{-1}$ can be implemented by solving the localized biharmonic problems on each subdomain in parallel. Hence, a preconditioner for $S_{h}^{-1}$ needs to be constructed.}
\end{remark}
\par
%
\section{A BDDC Preconditioner}\label{Sec:Precond_B2}
%
In this section a preconditioner for the Schur complement $S_h$ is constructed by the BDDC methodology.
\par
Let $V_{h,C,j}, 1\leq j\leq J$ be the restriction of $V_{h,C}$ on the subdomain $\D_j$. We define $\ExtH_{j}$, the space of local discrete biharmonic functions, by
\begin{equation*}\label{eq:HjDef}
  \ExtH_{j} = \left\{ v \in V_{h,C,j} : a_{h,j}(v,w) = 0  \quad  \forall\, w \in V_{h,C}(\D_{j})\right\},
\end{equation*}
where
$V_{h,C}(\D_{j})$ is the subspace of $V_{h,C,j}$ whose members vanish up to order $1$ on $\p\D_j$.
The space $\Hcc$ is then defined by gluing the spaces $\ExtH_j$ together at the cross points such that
\begin{equation*}\label{eq:HccDef}
  \Hcc = \left\{ v \in \L2 : v\big|_{\D_j} \in \ExtH_j \text{ and } v \text{ has continuous dofs at subdomain corners} \right\}.\notag
\end{equation*}
We equip $\Hcc$ with the bilinear form:
\begin{equation*}\label{eq:BDDC_BLform_Def}
\ahC(v,w) = \sum_{1 \leq j \leq J} a_{h,j}(v_j, w_j) \qquad \forall \,\, v, w \in \Hcc,
\end{equation*}
 where $v_j = v\big|_{\D_j}$ and $w_j = w\big|_{\D_j}$.
\par
 Next we introduce a decomposition of $\Hcc$,
$$\Hcc = \Hcz \oplus \Hc $$
where
\begin{align*}
  \Hcz &= \left\{ v \in \Hcc: \text{ the dofs of } v \text{ vanish at the corners of the subdomains } \D_1, \ldots, \D_J \right\}, \\
  \Hc &= \left\{ v \in \Hcc: a_h^{C}(v,w) = 0 \quad \forall \, w \in \Hcz \right\}. 
\end{align*}
\par
 Let $\Hcz_j$ be the restriction of $\Hcz$ on $\D_j$. We then define SPD operators $S_0:\Hc \longrightarrow \mathcal{H}'_{0}$ and $S_j:\Hcz_j\longrightarrow \Hcz_j'$  by
\begin{equation*}
\langle S_{0} v, w \rangle = \ahC(v,w) \quad \forall \, v,w \in \Hc \quad\text{and} \quad
\langle S_{j}v,w \rangle = a_{h,j}(v,w) \quad \forall \, v,w \in \Hcz_j.
\end{equation*}
\par
 Now the BDDC preconditioner $\PCS$ for $S_h$ is given by
\begin{equation*}\label{eq:PCSDef}
   \PCS= \left(\AvgP I_0 \right) S_{0}^{-1}\left(\AvgP I_0 \right)^{t} + \sum_{j=1}^{J} \left(\AvgP\bbE_j \right) S_{j}^{-1}\left(\AvgP\bbE_j \right)^{t},
\end{equation*}
where $I_0:\Hc \rightarrow \Hcc$ is the natural injection, $\bbE_j:\Hcz_j \rightarrow \Hcc$ is the trivial extension, 
and $\AvgP:\Hcc \longrightarrow V_{h,C}$ is a projection defined by averaging such that for all $v\in \Hcc, \, \AvgP v$ is continuous on $\Gamma$ up to order 1.
\par
\begin{remark}
 A preconditioner $B:\Vh{'}\longrightarrow \Vh$ for $A_h$ can then be constructed as follows:
\begin{equation*}\label{eq:B2Def}
  B = I_{D} A_{h,D}^{-1}I_{D}^{t} + I_{h,C,\D\backslash\Gamma} A_{h,C,\D\backslash\Gamma}^{-1} I_{h,C,\D\backslash\Gamma}^t + I_{\Gamma} \PCS I_{\Gamma}^t, 
\end{equation*}
where $I_D : V_{h,D} \rightarrow V_h, I_{h,C,\D\backslash\Gamma} : V_{h,C}(\D\backslash\Gamma) \rightarrow V_h$, and $I_{\Gamma} : V_{h,C}(\Gamma) \rightarrow V_{h}$ are natural injections.
\end{remark}

\section{Condition Number Estimates}\label{Sec:CondNoEst}
%
In this section we present the condition number estimates of $\PCS S_{h}$. Let us begin by noting that
\begin{equation*}\label{eq:ASDecomp}
  \VcG =  \AvgP I_{0} \Hc + \sum_{j=1}^{J} \AvgP \bbE_j \Hcz_j.
\end{equation*}
Then it follows from the theory of additive Schwarz preconditioners (see for example ~\cite{SBG:1996:DDM,TW:2005:DDM,Mathew:2008:DDM, BS:2008:FEM}) that the eigenvalues of $\PCS S_h$ are positive, and the extreme eigenvalues of $\PCS S_h$ are characteristic by the following formulas
\begin{align}
  \lambda_{\min}\left(\PCS S_h\right)
  &= \min_{\substack{\strut v \in \VcG \\ v \neq 0} }
  \frac{\langle S_h v, v \rangle}
  {
  \d\min_{\substack{\strut
    v = \AvgP I_0 v_0 + \sum_{j=1}^{J} \AvgP \bbE_{j} \vz_j \\
    v_0 \in \Hc, \vz_j \in \Hcz_j}}
  \Big(\langle S_0 v_0, v_0 \rangle +
  \sum_{j=1}^{J} \langle S_j  \vz_j, \vz_j  \rangle
  \Big)},
 \nonumber \\
  \lambda_{\max}\left(\PCS S_h\right)
   &= \max_{\substack{\strut v \in \VcG \\ v \neq 0}}
  \frac{\langle S_h v, v \rangle}
  {\d\min_{\substack{\strut
    v = \AvgP I_0 v_0 + \sum_{j=1}^{J} \AvgP \bbE_{j} \vz_j \\
    v_0 \in \Hc, \vz_j \in \Hcz_j}}
  \Big(\langle S_0 v_0, v_0 \rangle +
  \sum_{j=1}^{J} \langle S_j  \vz_j, \vz_j  \rangle
  \Big)}, \nonumber
\end{align}
from which we can establish a lower bound for the minimum eigenvalue of $\PCS S_h$, an upper bound for the maximum eigenvalue of $\PCS S_h$, and then an estimate on the condition number of $\PCS S_h$.
\begin{theorem}
 It holds that $\ds \lambda_{\min}(\PCS S_h) \geq 1$ and $\ds \lambda_{\max}(\PCS S_h) \leq (1+\ln(H/h))^2 / C$, which imply
  $$\kappa (\PCS S_h) = \frac{\lambda_{\min}(\PCS S_h)}{\lambda_{\max}(\PCS S_h)} \leq C (1+\ln(H/h))^2, $$
 where the positive constant $C$ is independent of $h, H$, and $J$.
\end{theorem}
%
\section{Numerical Results}\label{Sec:NumericalResults}
\noindent
In this section we present some numerical results to illustrate the performance of the preconditioners $\PCS$ and $B$. We consider our model problem \eqref{eq:model-weakform} on the unit square $(0,1)\times (0,1)$. By taking the penalty parameter $\eta$ in $a_h (\cdot, \cdot)$ and $a_{h, j}(\cdot, \cdot)$ to be 5, we compute the maximum eigenvalue, the minimum eigenvalue, and the condition number of the systems $\PCS S_h$ and $B A_h$ for different values of $H$ and $h$.
\par
The eigenvalues and condition numbers of $\PCS S_h$ and $B A_h$ for 16 subdomains are presented in Tables \ref{Table:PreSh16} and \ref{Table:PreAh16}, respectively. They confirm our theoretical estimates. In addition, the corresponding condition numbers of $A_h$ are provided in Table \ref{Table:PreAh16}. 
\par
Moreover, to illustrate the practical performance of the preconditioner, we present in Table \ref{Table:Niter_Ah16} the number of iterations required to reduce the relative residual error by a factor of $10^{-6}$ for the preconditioned system and the un-preconditioned system, from which we can observe the dramatic improvement in efficiency due to the preconditioner, especially as $h$ gets smaller.
\par
\begin{table}[h]
\centering
\caption{Eigenvalues and condition numbers of $\PCS S_h$ for $H=1/4$ ( J = 16 subdomains )}
\label{Table:PreSh16}
\begin{tabular}{l|c|c|c}
   \hspace{1.2cm} & $\quad \lambda_{\max}(\PCS S_h) \quad$ & $\quad \lambda_{\min}(\PCS S_h) \quad$ & $\quad \kappa (\PCS S_h) \quad$ \\[3pt]
 \hline \noalign{\smallskip}
 $h$=1/8 & 3.6073 & 1.0000 & 3.6073 \\[3pt]
 $h$=1/12 & 2.9197 & 1.0000 & 2.9197 \\[3pt]
 $h$=1/16 & 3.0908 & 1.0000 & 3.0908 \\[3pt]
 $h$=1/20 & 3.2756 & 1.0000 & 3.2756 \\[3pt]
 $h$=1/24 & 3.4535 & 1.0000 & 3.4535
\end{tabular}
\end{table}
\begin{table}[h]
\centering
\caption{Eigenvalues and condition numbers of $B A_h$, and condition numbers of $A_h$ for $H=1/4$ ( J = 16 subdomains )}
\label{Table:PreAh16}
\begin{tabular}{l|c|c|c|c}
   \hspace{1.2cm} & $\quad \lambda_{\max} (B A_h) \quad$ & $\quad \lambda_{\min}(B A_h) \quad$ & $\quad \kappa (B A_h) \quad$ & $\quad \kappa(A_h) \quad$\\[3pt]
 \hline  \noalign{\smallskip}
 $h$=1/8 & 4.0705 & 0.2148 & 18.9490 & 1.1064e+03 \\[3pt]
 $h$=1/12 & 3.4107 & 0.2507 & 13.6054 & 1.3426e+04 \\[3pt]
 $h$=1/16 & 3.4866 & 0.2578 & 13.5244 & 6.1689e+04
 \\[3pt]
 $h$=1/20 & 3.5947 & 0.2590 & 13.8787 & 1.8215e+05 \\[3pt]
 $h$=1/24 & 3.7123 & 0.2593 & 14.3181 & 4.2288e+05
\end{tabular}
\end{table}
\begin{table}[h]
\centering
\caption{Number of iterations for reducing the relative residual error by a factor of $10^{-6}$ for $H=1/4$ ( J = 16 subdomains )}
\label{Table:Niter_Ah16}
\begin{tabular}{l|c|c}
   \hspace{1.2cm} & $\quad Niter (A_h x = b) \quad$ & $\quad Niter (B A_h x = B b) \quad$ \\[3pt]
 \hline  \noalign{\smallskip}
 $h$=1/8 & 95 & 27 \\[3pt]
 $h$=1/12 & 235 & 23 \\[3pt]
 $h$=1/16 & 434 & 23 \\[3pt]
 $h$=1/20 & 704 & 23 \\[3pt]
 $h$=1/24 & 1026 & 23
\end{tabular}
\end{table}

\section*{Acknowledgements}
 The work of the first and third authors was supported in part by the National Science Foundation under Grant No. DMS-16-20273.

\end{document}